\documentclass{article}
\usepackage{amssymb,amsmath,amsthm}
\usepackage{graphicx,graphics,epsfig}
\newtheorem{theorem}{Theorem}
\newtheorem{definition}{Definition}

\newtheorem{proposition}{Proposition}

\newtheorem{remark}{Remark}
\newtheorem{example}{Example}
\date{}
\numberwithin{equation}{section} \numberwithin{theorem}{section}
\numberwithin{lemma}{section} \numberwithin{corollary}{section}
\numberwithin{remark}{section} \numberwithin{proposition}{section}
\numberwithin{definition}{section} \numberwithin{example}{section}

\newcommand{\vs}{\vskip}
\begin{document}

\title{Soft elementary compact in soft elementary topology}

  \vskip 0.5cm
\author{Mahmoud Bousselsal\\
Laboratory of Non Linear PDE  and History of Mathematics,\\
 ENS, 16050 Kouba, Algiers, Algeria.\\
Abderachid Saadi\\
Laboratory of Non Linear PDE  and History of Mathematics,\\
 ENS, 16050 Kouba, Algiers, Algeria.\\
Department of Mathematics,\\ Mohamed Boudiaf University of Msila,
B.P 166 Ichbelia, Msila, Algeria.\\} \maketitle

 \vs 0.5cm
\begin{abstract}
The notion of soft topology was introduced very recently, built up
on soft elementary intersection and union. In this paper,  Based on
this approach, we introduce the notion of soft elementary compact
sets and spaces. Also, we investigate their properties. To that end
we prove the soft elementary version of Baire theorem.
\end{abstract}

2010 AMS Classification: 54A05, 54A10, 54D30, 06D72.

 \vs 0.5cm


\section{Introduction}\label{1}

The theory of soft set is a new mathematical tool dealing with
uncertainties, it was introduced by Molodtsov \cite{[MOLODTSOV]} in
1999. This theory has several  applications in economy, medical
sciences,  social sciences, engineering,  see \cite{[Chen],[MAJI2]}.

After Molodtsov. Maji et al. defined some operators for soft sets
\cite{[MAJI1]}.  Later; those notions were improved  in \cite{[Ali],
[Sezgin], [Yang]}.

In \cite{[Shabir]}, Shabir and Naz defined the soft topology, based
on the intersection and union of soft sets as in \cite{[Sezgin]}.
This notion of soft topology is used in \cite{[Aras], [Aygunoglu],
[Cagman], [Georgiou1], [Georgiou2], [Min], [Zorlutuna]}. Soft
concepts of element, interior, closure, separtion were also defiend
in \cite{[Shabir]}.

Hazra et al. \cite{[Hazra]} gave two approaches of soft topology,
and deduce some properties and results.

The notion of soft element defined by Das and Samanta
\cite{[Das-Samanta3]}, is compatible with the definition of soft
subsets, and different from that one defined in  \cite{[Shabir]}. It
is used to define soft real number, soft complex and soft metric
spaces , see \cite{[Das-Samanta1], [Das-Samanta2], [Das-Samanta3]}.
From this notion of soft element, Das and Samanta introduced other
operations such  that elementary union; intersection, and
complement.

The compactness was investigated by Aygunoglu and Aygun
\cite{[Aygunoglu]}, Bayramov and Gunduz \cite{[Bayramov-Gunduz]},
and  Zorlutuna et al \cite{[Zorlutuna]}, which are based at the
definition given by Shabir in \cite{[Shabir]}.

In \cite{[Chiney-Samanta]}, Chiney and Samanta introduced a new
definition of soft topology based on the elementary intersection,
elementary union, elementary complement. This soft topology is
different from these defined in \cite{[Shabir], [Hazra]}. The
necessary soft topological concept tools needed for continuity are
defined (soft interior element and set, limiting soft, soft closure
set, soft neighborhood, soft base, soft continuous function and soft
separation axioms).

In this paper, we adopt this new definition, and for spare a
probably paradox, we denote by soft elementary topology (or soft
e-topology) for the soft topology which is defined in
\cite{[Chiney-Samanta]}.

The paper comprised five sections. In section 2, we introduce  some
well-known results in soft set theory, which are needed in the
paper.  In section 3, we introduce the definition  of soft
elementary  sub-topology,  and some other properties and results. In
the section 4, we state the definition of soft elementary
quasi-compact space, soft elementary compact space, soft elementary
compact set and their proprieties. The last section is devoted to
the proof of the soft version of Baire theorem.

\vs 0.5cm


\section{Preliminaries}\label{2}
This section contains definitions and properties of soft sets, and
some related notions that will be needed throughout this paper. Let
$X$ be an initial universe set, $E$ be the set of parameters, and
$A$ be a nonempty subset of $E$.
\begin{definition}\textsc{\cite{[MOLODTSOV]}}
A pair $(F,A)$ is called a soft set over $X$, if and only if $F$ is
a mapping of $A$ into $\mathcal{P}(X)$.
\end{definition}

\begin{definition}\textsc{\cite{[Ali], [MAJI1], [Sezgin]}}
Let $(F,A)$ and $(G,A)$ be two soft sets over $X$.
\begin{description}
    \item[i)] $(F,A)$ is called a soft subset of $(G,A)$, $($ i.e. $(F,A)\widetilde{\subseteq}(G,A)),$ if $F(\alpha)$ is a subset
of $G(\alpha)$, for all $\alpha\in A$;
    \item[ii)] $(F,A)$ and $(G,A)$ are called soft equal if $(F,A)$ is a
soft subset of $(G,A)$ and $(G,A)$ is a soft subset of $(F,A)$;
    \item[iii)] the complement or relative complement of a soft set
$(F,A)$ is denoted by $(F,A)^C$ and is defined by $(F,A)^C=(F^C,A)$,
 where $F^C(\alpha) =C_X^{F(\alpha)}$, for all $\alpha\in A$.
    \item[iv)] union of $(F,A)$ and $(G,A)$ is the soft set $(H,A),$ defined by
$H(\alpha) = F(\alpha)\cup G(\alpha)$, for all $\alpha\in A,$ and
    denoted by $(F,A)\widetilde{\cup}(G,A)$;
    \item[v)] intersection of  $(F,A)$ and $(G,A)$ is the soft set $(H,A),$ defined by
$H(\alpha) = F(\alpha)\cap G(\alpha)$, for all $\alpha\in A,$ and
    denoted by $(F,A)\widetilde{\cap}(G,A)$.
\end{description}
\end{definition}

\begin{example}\label{ex2.1}
Let $X=\{x,y,z\}, A=\{\alpha, \beta\}, (F,A)$ and $(G,A)$ such that
$F(\alpha)=\{x,y\}, F(\beta)=\{x,z\}, G(\alpha)=\{y,z\},
G(\beta)=\{x\}$. Then; $(F, A)\widetilde{\cup}(G,A)=(H,A)$ such that
$H(\alpha)=X, H(\beta)=\{x,z\}$, $(F, A)\widetilde{\cap}
(G,A)=(J,A)$ such that $J(\alpha)=\{y\}, H(\beta)=\{x\}$, $(F,
A)^c=(F^C,A)$ such that $F^c(\alpha)=\{z\}, F^C(\beta)=\{y\}$.
\end{example}

\begin{definition}\textsc{\cite{[Shabir]}}
Let $Y$ be a subset of $X$. We denote by $(\widetilde{Y},A)$ the
soft set $(F,A)$ such that $F(\alpha)=Y,$ for all $\alpha\in A.$ If
$Y=\emptyset, (F,A)$ is called null soft, set denoted by
$(\widetilde{\Phi},A)$.  If $Y = X, (F,A)$ is called absolute soft
set, denoted by $(\widetilde{X},A)$. $S(\widetilde{X})$ is the
collection of $(\widetilde{\Phi},A)$, and the soft sets $(F,A)$ such
that $(F,A)(\alpha)\neq\emptyset$, for all $\alpha\in A$.
\end{definition}

It is obvious that if $Y\neq\emptyset,$ then $(\widetilde{Y},A)\in
S(\widetilde{X}),$ and $(\widetilde{Y},A)\neq (\widetilde{\Phi},A)$.\\
We denote by $S(\widetilde{Y})$ the collection of soft subsets
$(F,A)$ of $(\widetilde{Y},A)$ such that:
$(F,A)=(\widetilde{\Phi},A)$ or $(F,A)(\alpha)\neq\emptyset$ for all
$\alpha\in A.$

\begin{remark}\label{rem2.1}
Let $Y$ be a nonempty subset of $X,$ and $(F,A)\in
S(\widetilde{X})$. If $(F,A)\widetilde{\cap} (\widetilde{Y},A)\in
S(\widetilde{X})$ then $(F,A)\widetilde{\cap} (\widetilde{Y},A)\in
S(\widetilde{Y})$. Indeed, $(F,A)\widetilde{\cap}
(\widetilde{Y},A)\widetilde{\subset} (\widetilde{Y},A),$ and
$(F,A)=(\widetilde{\Phi},A)$ or $F(\alpha)\neq\emptyset$ for all
$\alpha\in A.$ Hence $(F,A)\widetilde{\cap} (\widetilde{Y},A)\in
S(\widetilde{Y}).$
\end{remark}

\begin{definition}\textsc{\cite{[Das-Samanta2]}}
\begin{description}
    \item[i)] A soft element of $(\widetilde{X}, A)$ is a function $\widetilde{x},$ defined on $A$ to the set $X.$ A soft element
    $\widetilde{x}$ of $(\widetilde{X}, A)$ is said to belong to a soft set $(F,A)$ over $X$, which is denoted by
$\widetilde{x}\widetilde{\in}(F,A)$, if $\widetilde{x}(\alpha)\in
F(\alpha)$,  for all $\alpha\in A$. Therefore, if $(F,A)$ is such
that $F(\alpha)\neq\emptyset$, for all $\alpha \in A$, we have
$F(\alpha) = \{\widetilde{x}(\alpha):
\widetilde{x}\widetilde{\in}(F,A)\}$, for all $\alpha\in A$;
    \item[ii)] the collection of all soft elements of a soft set
$(F,A)$ is denoted by $SE(F,A)$;
    \item[iii)] for a collection $B$ of soft elements of $(\widetilde{X},
    A)$, we denote by $SS(\mathcal{B})$ the soft set $(F,A)$ such that: $F(\alpha)=\{\widetilde{x}(\alpha), \widetilde{x}\in \mathcal{B}\}.$
\end{description}
\end{definition}

\begin{definition}\textsc{\cite{[Das-Samanta2]}}
For any two soft sets $(F,A), (G,A) \in S(\widetilde{X})$. Then
\begin{description}
    \item[i)] elementary union of $(F,A)$ and $(G,A)$ is denoted by $(F,A)\Cup(G,A),$ and defined by
    $(F,A)\Cup(G,A) = SS(SE(F,A)\cup SE(G,A))$;
    \item[ii)] elementary intersection of $(F,A)$ and $(G,A)$ is denoted by $(F,A)\Cap(G,A),$ and defined by $(F,A)\Cap(G,A) = SS(SE(F,A)\cap
    SE(G,A))$;
    \item[iii)] elementary
complement of $(F,A)$ is denoted by $(F,A)^\mathbb{C},$ and defined
by $(F,A)^\mathbb{C} = SS(\mathcal{B})$, where $\mathcal{B} =
\{\widetilde{x}\widetilde{\in}(\widetilde{X},A): \widetilde{x}\in
(F,A)^C\}$.
\end{description}
\end{definition}

\begin{example}\label{ex2.2}
Let $X=\{x,y,z\}, A=\{\alpha, \beta\}.$ Then; $(\widetilde{X},A)=
SS(\{\widetilde{x}_1, \widetilde{x}_2, \widetilde{x}_3,
\widetilde{x}_4, \widetilde{x}_5, \widetilde{x}_6,$\\
$\widetilde{x}_7, \widetilde{x}_8, \widetilde{x}_9\})$, where:
$\widetilde{x}_1(\alpha)=\widetilde{x}_1(\beta)=x,
\widetilde{x}_2(\alpha)=\widetilde{x}_2(\beta)=y,
\widetilde{x}_3(\alpha)=\widetilde{x}_3(\beta)=z,
\widetilde{x}_4(\alpha)=x, \widetilde{x}_4(\beta)=y,
\widetilde{x}_5(\alpha)=x, \widetilde{x}_5(\beta)=z,
\widetilde{x}_6(\alpha)=y, \widetilde{x}_6(\beta)=x,
\widetilde{x}_7(\alpha)=y, \widetilde{x}_7(\beta)=z,
\widetilde{x}_8(\alpha)=z, \widetilde{x}_8(\beta)=x,
\widetilde{x}_9(\alpha)=z, \widetilde{x}_4(\beta)=y.$ Let $(F,A),
(G,A)\in S(\widetilde{X})$ such that $F(\alpha)=\{x,z\},
F(\beta)=\{y,z\}, G(\alpha)=\{x,y\}, G(\beta)=\{x\}$. Then,
$(F,A)=SS(\{\widetilde{x}_3, \widetilde{x}_4, \widetilde{x}_5,
\widetilde{x}_9\}),$\\ $(G,A)= SS(\{\widetilde{x}_1,
\widetilde{x}_6\})$. Hence, $(F, A)\Cup(G,A)= SS(\{\widetilde{x}_1,
\widetilde{x}_3, \widetilde{x}_4, \widetilde{x}_5, \widetilde{x}_6,
\widetilde{x}_9\}), (F, A)\Cap(G,A)= (\widetilde{\Phi},A)$.  We have
$(F, A)\Cup(G,A)=(F,A)\widetilde{\cup}(G,A)=(\widetilde{X},A)$,  but
$(F, A)\Cap(G,A)\neq(F,A)\widetilde{\cap}(G,A)$ since
$(F,A)\widetilde{\cap}(G,A)= (H,A),$ such that $H(\alpha)=\{x\},
H(\beta)=\emptyset$. Now, $(F,A)^C=(F,A)^\mathbb{C}=
SS(\{\widetilde{x}_1, \widetilde{x}_2, \widetilde{x}_6,
\widetilde{x}_8\}).$
\end{example}

\noindent\textbf{Notation:}  Let $Y$ be a nonempty subset of $X$ and
let $(Z,A)\in S(\widetilde{Y})$. We denote by $(Z,A)_Y^C$ the soft
set $(W,A)$ over $Y,$ where $W(\alpha)=Y\setminus Z(\alpha)$ for all
$\alpha\in A.$ Also,we denote by $(Z,A)_Y^\mathbb{C}$ the soft set
of the soft elements $\widetilde{x}$ such that $\widetilde{x}
\widetilde{\in} (Z,A)_Y^C.$

\begin{definition}\textsc{\cite{[Chiney-Samanta]}}\label{def2.6}
Let $\tau$ be a collection of soft sets of $S(\widetilde{X})$. Then
$\tau$ is called a soft e-topology on $(\widetilde{X},A)$ if the
following conditions hold
\begin{description}
    \item[i)] $(\widetilde{\Phi},A)$ and $(\widetilde{X},A)$ belong
    to   $\tau$;
    \item[ii)] the elementary union of any number of soft sets in $\tau$ belongs to
    $\tau$;
    \item[iii)] the elementary intersection of two soft sets in $\tau$ belongs to $\tau$.
\end{description}
The triplet $(\widetilde{X},\tau,A)$ is called a soft e-topological
space. A member of $\tau$ is called soft e-open sets in
$(\widetilde{X},\tau,A)$.
\end{definition}

\begin{definition}\textsc{\cite{[Chiney-Samanta]}}
A soft set $(F,A)\in S(\widetilde{X})$ is called a soft e-closed set
in $(\widetilde{X},\tau,A)$ if its relative complement $(F,A)^C$
belongs to $S(\widetilde{X})$ and $(F,A)^\mathbb{C}$ belongs to
$\tau $.
\end{definition}

\begin{proposition}\textsc{\cite{[Chiney-Samanta]}}
\begin{description}
    \item[i)] $(\widetilde{\Phi},A)$ and $(\widetilde{X},A)$ are soft e-closed soft
sets in $(\widetilde{X},\tau,A)$,
    \item[ii)] arbitrary elementary intersection of soft e-closed sets is soft e-closed.
\end{description}
\end{proposition}

\begin{remark}\textsc{\cite{[Chiney-Samanta]}}
In general, the elementary union of two soft e-closed sets is not
soft e-closed.
\end{remark}

\begin{definition}\textsc{\cite{[Chiney-Samanta]}}
Let $(\widetilde{X},\tau,A)$  be a soft e-topological space and
$(F,A) \in S(\widetilde{X})$. Then the soft e-closure of $(F,A)$,
denoted by $\overline{(F,A)}$ is defined as the elementary
intersection of all soft e-closed super sets of $(F,A)$.
\end{definition}

\begin{definition}\textsc{\cite{[Chiney-Samanta]}}
Let $(\widetilde{X},\tau,A)$ be a soft e-topological space. A soft
element $\widetilde{x}\in(\widetilde{X},A)$ is called a limiting
soft element of a soft set $(F,A)\in S(\widetilde{X}),$ if for all
$(G,A)\in \tau,$ and for any $\alpha\in A$,
$\widetilde{x}(\alpha)\in G(\alpha).$ This implies that $F(\alpha)
\cap [G(\alpha)\setminus \widetilde{x}(\alpha)] \neq\emptyset$.
\end{definition}

\begin{definition}\textsc{\cite{[Chiney-Samanta]}}
Let $(\widetilde{X},\tau,A)$  be a soft e-topological space and
$(F,A)\in S(\widetilde{X})$. A soft element
$\widetilde{x}\widetilde{\in}(F,A)$ is called an interior soft
element of $(F,A),$ if there exists $(G,A) \in \tau$ such that
$\widetilde{x}\widetilde{\in}(G,A)\widetilde{\subset}(F,A)$. The
interior of a soft set $(F,A)$ denoted by $Int(F,A)$, is defined by
$Int(F,A) = \{\widetilde{x}\widetilde{\in}(F,A):
\widetilde{x}\widetilde{\in}(G,A)\widetilde{\subset} (F,A) \,for\,
some (G,A)\in\tau\}.  SS[Int(F,A)]$ is called soft interior of
$(F,A)$ and denoted by $\overset{\circ}{\overbrace{(F,A)}}$.
\end{definition}

\begin{definition}\textsc{\cite{[Chiney-Samanta]}}
Let $(\widetilde{X},\tau,A)$ be a soft e-topological space. Then
$(\widetilde{\Phi},A)\neq(F,A)\in S(\widetilde{X})$ is a soft
neighborhood (soft nbd, for short) of the soft element
$\widetilde{x}$ if there exists a soft set $(G,A)\in\tau,$ such that
$\widetilde{x}\widetilde{\in}(G,A)\widetilde{\subset}(F,A)$.
\end{definition}

\begin{example}\label{ex2.3}
Let $X=\{a,b,c,d\}, A=\{\alpha,\beta\},$ and
$\tau=\{(\widetilde{\Phi},A), (\widetilde{X},A), (F_1,A),
(F_2,A),$\\$ (F_3,A), (F_4,A)\},$ such that $F_1(\alpha)=\{a\},
F_1(\beta)=\{b\}, F_2(\alpha)=\{b,c\}, F_2(\beta)=\{c,d\},$\\
$F_3(\alpha)=\{a,b,c\}, F_3(\beta)=\{b,c,d\}, F_4(\alpha)=X,
F_4(\beta)=\{b,c,d\}$. Then; $(\widetilde{X},\tau,A)$ is a soft
e-topological space, and the collection of e-closed  sets is
$\{(\widetilde{\Phi},A), (\widetilde{X},A), (F_1,A)^\mathbb{C},
(F_2,A)^\mathbb{C},$\\ $(F_3,A)^\mathbb{C},\}.$ $(F_4,A)^C$ is not a
soft e-closed set, since $(F_4,A)^C\in S(\widetilde{X})$.\\ Now, let
$(F,A), (G,A)$ be such that $F(\alpha)=F(\beta)=\{c\},
G(\alpha)=\{a,b\}, G(\beta)=\{b\}.$ Then;
$\overline{(F,A)}=(F_1,A)^\mathbb{C},
\overset{0}{\overbrace{(G,A)}}=(F_1,A),$ and $(G,A)$ is a soft  nbd
of $(\widetilde{x})$ such that $\widetilde{x}(\alpha)=a,
\widetilde{x}(\beta)=b$.
\end{example}

\begin{definition}\textsc{\cite{[Chiney-Samanta]}}
Let $(\widetilde{X},\tau,A)$ be a soft e-topological space. Let
$\widetilde{x}, \widetilde{y} \in (\widetilde{X},A)$ such that
$\widetilde{x}(\alpha)\neq \widetilde{y}(\alpha),$ for all
$\alpha\in A$. Then, if  there exist $(F,A), (G,A) \in \tau,$ such
that $\widetilde{x}\in (F,A), \widetilde{y}\in (G,A)$ and
$(F,A)\widetilde{\cap}(G,A) = (\widetilde{\Phi},A)$, then
$(\widetilde{X},\tau,A)$  is called a soft $e-T_2$ space,  or soft
e-Hausdorff space.
\end{definition}

\begin{definition}\textsc{\cite{[Chiney-Samanta]}}
A soft e-topological space $(\widetilde{X},\tau,A)$  is called a
soft e-regular space if for any soft closed set $(F,A)$ and for any
soft element $\widetilde{x},$ such that $\widetilde{x}(\alpha)\notin
(F,A)(\alpha)$, for all $\alpha\in A$, there exist $(G,A), (H,A) \in
\tau,$ such that $(F,A)\widetilde{\subset} (G,A),
\widetilde{x}\widetilde{\in}(H,A)$ and $(F,A)\Cap(G,A) =
(\widetilde{\Phi},A)$.\end{definition}

\begin{definition}\textsc{\cite{[Chiney-Samanta]}}
A soft e-topological space $(\widetilde{X},\tau,A)$  is called a
soft e-normal space if for any two soft closed sets $(F,A)$ and
$(G,A),$ such that $(F,A)\widetilde{\cap}(G,A) =
(\widetilde{\Phi},A)$, there exist $(U,A), (V,A) \in \tau,$ such
that $(F,A)\widetilde{\subset}(U,A), (G,A)\widetilde{\subset}(V,A)$
and $(U,A)\Cap (V,A)=(\widetilde{\Phi},A)$.
\end{definition}

\begin{definition}\textsc{\cite{[Chiney-Samanta]}}
Let $X, Y$ be two non-empty sets and $\{f_\alpha: X\rightarrow Y,
\alpha\in A\}$  be a collection of functions. Then a function $f:
SE(\widetilde{X})\rightarrow SE(\widetilde{Y})$ defined by
$[f(\widetilde{x})](\alpha) = f_\alpha(x(\alpha)),$ for all
$\alpha\in A$ is called a soft function.
\end{definition}

\begin{definition}\textsc{\cite{[Chiney-Samanta]}}
Let  $f: SE(\widetilde{X})\rightarrow SE(\widetilde{Y})$ be a soft
function. Then
\begin{description}
    \item[i)] the image of a soft set $(F,A)$ over $X $ under the soft function $f$, denoted by
$f[(F,A)]$, is defined by $f[(F,A)] = SS\{f(SE(F,A))\},$ i.e.
$f[(F,A)](\alpha) = f_\alpha(F(\alpha))$, for all $\alpha\in A$.
    \item[ii)] the inverse image of a soft set $(G,A)$ over $Y$ under the soft function $f$, denoted by $f^{-1}[(G,A)]$, is defined by $f^{-1}[(G,A)]=
    SS\{f^{-1}(SE(G,A))\}$, i.e. $f^{-1}[(G,A)](\alpha) = f^{-1}_\alpha (G(\alpha))$, for all $\alpha\in A$.
\end{description}
\end{definition}

\begin{definition}\textsc{\cite{[Chiney-Samanta]}}
Let $(\widetilde{X},\tau,A), (\widetilde{Y},\sigma,A)$ be two soft
e-topological spaces and\\ $f: SE(\widetilde{X}) \rightarrow
SE(\widetilde{Y})$ be a soft function. $f:
(\widetilde{X},\tau,A)\rightarrow (\widetilde{Y},\sigma,A)$ is
called soft e-continuous at $\widetilde{x}_0 \in(\widetilde{X},A)$,
if for every $(V,A)\in \sigma$ such that
$f(\widetilde{x}_0)\widetilde{\in}(V,A)$, there exists $(U,A)
\in\tau$ such that $\widetilde{x}_0\widetilde{\in}(U,A)$ and
$f(U,A)\widetilde{\subseteq}(V,A)$. $f$ is called soft e-continuous
on $(\widetilde{X},\tau,A),$ if it is soft e-continuous at each soft
element $\widetilde{x}_0 \in(\widetilde{X},A)$.
\end{definition}

\begin{proposition}\textsc{\cite{[Chiney-Samanta]}}
Let $(\widetilde{X},\tau,A), (\widetilde{Y},\sigma,A)$ be two soft
e-topological spaces and\\ $f: SE(\widetilde{X}) \rightarrow
SE(\widetilde{Y})$ be a soft function. Then, $f$ is  soft
e-continuous on $(\widetilde{X},\tau,A),$ if and only if for any
soft e-open  set $(U,A)\in S(\widetilde{Y})$ in
$(\widetilde{Y},\sigma,A), f^{-1}(U,A)$ is a soft e-open set in
$(\widetilde{X},\tau,A)$.
\end{proposition}

\vs 0.5cm


\section{Soft elementary topology}\label{3}

In this section, based on the definition \ref{def2.6}, we introduce
the notion of soft sub-e-topological space, and we investigate their
properties. First, we start with the following key result.

\begin{theorem}
Let $(\widetilde{X},\tau, A)$ be a soft e-topological space such
that for all $(O_1,A), (O_2,A)\in \tau$ we have,
$(O_1,A)\widetilde{\cap} (O_2,A)\in S(\widetilde{X}).$ Let $Y$ be a
nonempty subset of $X$ such that for all $(O,A)\in \tau$ we have
$(O,A)\widetilde{\cap} (\widetilde{Y},A)\in S(\widetilde{X}).$ Then,
it holds that the collection
$\tau_Y=\{(O_Y,A)=(O,A)\Cap(\widetilde{Y},A), (O,A)\in \tau\}$
define a soft e-topology for $(\widetilde{Y},A).$
\end{theorem}

\noindent\emph{Proof.} Let $(\widetilde{X},\tau, A)$ be a soft
e-topological space and let $Y$ be a nonempty subset of $X$. We show
that $\tau_Y$ satisfies the conditions of soft e-topology. On the
one hand, since
$(\widetilde{\Phi},A)=(\widetilde{\Phi},A)\Cap(\widetilde{Y},A)$ and
$(\widetilde{Y},A)=(\widetilde{X},A)\Cap(\widetilde{Y},A)$.
    Then $(\widetilde{\Phi},A)$ and $(\widetilde{Y},A)$ belong  to
    $\tau_Y$.
On the other hand, suppose that
 $\{(O^i_Y,A)=(O^i,A)\Cap(\widetilde{Y},A), i\in
    I\}$ is a family of soft sets in $\tau_Y$, Then  $\{(O^i,A), i\in
    I\}$ is a family of soft sets in $\tau,$ and $(O,A)=\underset{i\in I}{\Cup}(O^i,A)\in
    \tau, $ it follows that\\
    $\begin{array}{lll}
     \underset{i\in I}{\Cup} (O^i_Y,A) & = & \underset{i\in I}{\Cup}[(O^i,A)\Cap(\widetilde{Y},A)] \\
                              & = & \underset{i\in I}{\widetilde{\cup}}[(O^i,A)\Cap(\widetilde{Y},A)] \\
                              & = & \underset{i\in I}{\widetilde{\cup}}[(O^i,A)\widetilde{\cap}(\widetilde{Y},A)] \\
                              & = & [\underset{i\in I}{\widetilde{\cup}}(O^i,A)]\widetilde{\cap}(\widetilde{Y},A) \\
                              & = & [\underset{i\in I}{\Cup}(O^i,A)]\widetilde{\cap}(\widetilde{Y},A) \\
                              & = & (O,A)\widetilde{\cap}(\widetilde{Y},A) \\
                              & = & (O,A)\Cap(\widetilde{Y},A) \in \tau_Y.
    \end{array}$\\
Finally, let $(O^1_Y,A)=(O^1,A)\Cap(\widetilde{Y},A),
    (O^2_Y,A)=(O^2,A)\Cap(\widetilde{Y},A)$ be two soft sets in $\tau_Y$. Then,  $(O^1,A), (O^2,A)$ are two soft sets in $\tau,$
     and  $(O,A)=(O^1,A)\Cap(O^2,A)\in\tau,$ it follows that\\
    $\begin{array}{lll}
      (O^1_Y,A)\Cap (O^2_Y,A) & = & [(O^1,A)\Cap(\widetilde{Y},A)]\Cap[(O^2,A)\Cap(\widetilde{Y},A)] \\
                              & = &
                              [(O^1,A)\widetilde{\cap}(\widetilde{Y},A)]\Cap[(O^2,A)\widetilde{\cap}(\widetilde{Y},A)].
\end{array}$\\
If $[(O^1,A)\widetilde{\cap}(\widetilde{Y},A)]
\Cap[(O^2,A)\widetilde{\cap}(\widetilde{Y},A)]=(\widetilde{\Phi},A)$
then $(O^1_Y,A)\Cap (O^2_Y,A)\in \tau_Y,$ else we obtain\\
    $\begin{array}{lll}
      (O^1_Y,A)\Cap (O^2_Y,A) & = & [(O^1,A)\widetilde{\cap}(\widetilde{Y},A)]\widetilde{\cap}[(O^2,A)\widetilde{\cap}(\widetilde{Y},A)]  \\
                              & = & [(O^1,A)\widetilde{\cap}(O^2,A)]\widetilde{\cap}(\widetilde{Y},A) \\
                              & = & (O,A)\Cap(\widetilde{Y},A) \in \tau_Y.
    \end{array}$\\
Thus, the collection $\tau_Y$ define a soft topology for
$(\widetilde{Y},A)$.

\noindent Next, we introduce the notion of sub-e-topological space,
based on the above theorem.

\begin{definition}
Let $(\widetilde{X},\tau, A)$ be a soft e-topological space such
that for all $(O_1,A), (O_2,A)\in \tau$ we have,
$(O_1,A)\widetilde{\cap} (O_2,A)\in S(\widetilde{X}).$ Let $Y$ be a
nonempty subset of $X$ such that for all $(O,A)\in \tau, $
$(O,A)\widetilde{\cap} (\widetilde{Y},A)\in S(\widetilde{X}),$ and
$\tau_Y=\{(O_Y,A)=(O,A)\Cap(\widetilde{Y},A), (O,A)\in \tau\}.$ The
triplet $(\widetilde{Y},\tau_Y, A)$ is called a soft
sub-e-topological space of $(\widetilde{X},\tau, A)$ and $\tau_Y$ is
called soft sub-e-topology of $\tau$. The members of $\tau_Y$ are
called soft $e_Y$-open sets in $(\widetilde{Y},\tau_Y,A)$.
\end{definition}

\begin{definition}
Let $(\widetilde{X},\tau, A)$ be a soft e-topological space such
that for all $(O_1,A), (O_2,A)\in \tau$ we have,
$(O_1,A)\widetilde{\cap} (O_2,A)\in S(\widetilde{X}),
(\widetilde{Y},\tau_Y, A)$ be a soft sub e-topological space of
$(\widetilde{X},\tau, A),$ and  $(Z,A)\in S(\widetilde{Y})$. The
soft set $(Z,A)$ is called a soft $e_Y$-closed set in
$(\widetilde{Y},\tau_Y,A),$ if $(Z,A)_Y^C \in S(\widetilde{Y})$ and
$(Z,A)_Y^\mathbb{C}\in \tau_Y.$
\end{definition}

\begin{example}\label{ex3.1}
Let $X=\{a,b,c,d\}, A=\{\alpha,\beta\},$ and
$\tau=\{(\widetilde{\Phi},A), (\widetilde{X},A), (F,A), (G,A),
(H,A)\},$ such that $F(\alpha)=\{a\}, F(\beta)=\{c,d\},
G(\alpha)=\{c,d\}, G(\beta)=\{a\}, H(\alpha)=\{a,c,d\},
H(\beta)=\{a,c,d\}$. Then, $(\widetilde{X},\tau, A)$ is a soft
e-topological space.\\ Let $Y=\{a,c\},$ then $\tau_Y=
\{(\widetilde{\Phi},A), (\widetilde{Y},A),(F_Y,A), (G_Y,A)\},$ is a
soft topology at $(\widetilde{Y},A),$ where $F_Y(\alpha)=\{a\},
F_Y(\beta)=\{c\}, G_Y(\alpha)=\{c\}, G_Y(\beta)=\{a\}$.
\end{example}

In the following proposition, we characterize the soft $e_Y-$closed
sets in a sub-e-topological space $(\widetilde{Y},\tau_Y, A)$.

\begin{proposition} Let $(\widetilde{X},\tau, A)$ be a soft e-topological space such that
for all $(O_1,A), (O_2,A)\in \tau$ we have, $(O_1,A)\widetilde{\cap}
(O_2,A)\in S(\widetilde{X}),$ and $(\widetilde{Y},\tau_Y, A)$ be a
soft sub-e-topological space of $(\widetilde{X},\tau, A)$.  If the
soft set $(Z,A)$ is a soft $e_Y$-closed set  in
    $(\widetilde{Y},\tau_Y,A),$ then there exists a soft e-closed
    set $(F,A)$ in $(\widetilde{X},\tau, A),$ such that $(Z,A)=(F,A)\Cap(\widetilde{Y},A).$
\end{proposition}

\noindent\emph{Proof.} Let $(Z,A)$ be a soft $e_Y$-closed set in
$(\widetilde{Y},\tau_Y,A).$ There are two cases to be considered
\begin{description}
    \item[Case (1)] If $(Z,A)=(\widetilde{\Phi},A),$ then
    $(F,A)=(\widetilde{\Phi},A).$
    \item[Case (2)] If $(Z,A)\neq(\widetilde{\Phi},A),$ then there exists
    $(O,A)\in \tau$ such that,
    $(Z,A)_Y^\mathbb{C}=(Z,A)_Y^C=(O,A)\Cap(\widetilde{Y},A).$ Hence, for all $\alpha\in A,$ we have $Y\setminus Z(\alpha)=Y\cap
    O(\alpha).$  Then, it follows that $Z(\alpha)=O^C(\alpha)\cap Y,$ for all $\alpha\in A$. Since, $(Z,A)\neq(\widetilde{\Phi},A),$ we get
    $(O,A)^C\in S(\widetilde{X}), O^C(\alpha)\neq\emptyset$ for all
    $\alpha\in A$ and $(Z,A)_Y^C(\alpha)=O(\alpha)^C\cap Y$. Putting $(F,A)=(O,A)^C,$ thus $(Z,A)_Y^C=(F,A)\widetilde{\cap}(\widetilde{Y},A)=(F,A)\Cap(\widetilde{Y},A).$
\end{description}

\begin{proposition}\label{lem3.1}
Let $(\widetilde{X},\tau,A)$  be a soft e-topological space, $(F,A)$
a soft subset of $(\widetilde{X},A)$, and
$\widetilde{x}\widetilde{\in} (\widetilde{X},A)$. If $\widetilde{x}$
is a soft limiting element of $(F,A),$ for all $(G,A)\in \tau$;
$\widetilde{x}\widetilde{\in}(G,A)$ implies that there exists
$\widetilde{y}\widetilde{\in} SE(\widetilde{X})$ such that
$\widetilde{y}\neq\widetilde{x}$ and $\widetilde{y}\widetilde{\in}
(F,A)\Cap(G,A)$.
\end{proposition}

\noindent\emph{Proof.} Let $(\widetilde{X},\tau,A)$  be a soft
topological space, $(F,A)$ a soft subset of $(\widetilde{X},A)$, and
$\widetilde{x}$ is a soft limiting element of $(F,A)$. Then,  for
any $(G,A)\in \tau$ and for any $\alpha\in A$,
$\widetilde{x}(\alpha)\in G(\alpha)$ implies that $F(\alpha) \cap
[G(\alpha) \setminus \{\widetilde{x}(\alpha)\}] \neq\emptyset$.
Hence, there exists $a_\alpha\in X$ such that,
$\widetilde{x}(\alpha)\notin F(\alpha)\cap G(\alpha)$. Let
$\widetilde{y}\in (\widetilde{X},A)$ such that
$\widetilde{y}(\alpha)= a_\alpha$ for all $\alpha\in A,$ then
$\widetilde{y}\neq\widetilde{x}$ and $\widetilde{y}\widetilde{\in}
(F,A)\Cap(G,A)$.

In the following proposition, we show that the sub-e-topological
space of an e-Hausdorff space is an e-Hausdorff space.

\begin{proposition} Let $(\widetilde{X},\tau, A)$ be a soft e-topological space such that
for all $(O_1,A), (O_2,A)\in \tau$ we have, $(O_1,A)\widetilde{\cap}
(O_2,A)\in S(\widetilde{X}),$ and $(\widetilde{Y},\tau_Y, A)$ be a
soft sub-e-topological space of $(\widetilde{X},\tau, A)$.
 If $(\widetilde{X},\tau, A)$ is a soft e-Hausdorff space then $(\widetilde{Y},\tau_Y, A)$
 is a soft e-Hausdorff space.
\end{proposition}

\noindent\emph{Proof.} Let  $\widetilde{x}, \widetilde{y} \in
(\widetilde{Y},A)$ such that $\widetilde{x}(\alpha)\neq
\widetilde{y}(\alpha)$ for all $\alpha\in A,$ then $\widetilde{x},
\widetilde{y} \in (\widetilde{X},A)$ and since $(\widetilde{X},\tau,
A)$ is a soft e-Hausdorff space there exist $(F,A), (G,A) \in \tau$
such that $\widetilde{x}\in (F,A), \widetilde{y}\in (G,A)$ and
$(F,A)\widetilde{\cap}(G,A) =
(\widetilde{\Phi},A)$.\\
$\widetilde{x}\in (F_Y,A)=(F,A)\Cap(\widetilde{Y},A),
\widetilde{y}\in (G_Y,A)=(G,A) \Cap(\widetilde{Y},A)$ and
$(F_Y,A)\widetilde{\cap}(G_Y,A) = (\widetilde{\Phi},A)$. Then
$(\widetilde{Y},\tau_Y, A)$
 is a soft e-Hausdorff space.

\vs 0.5cm


\section{Soft e-compact space and soft e-compact set}\label{4}
This section contains basic definitions and properties of soft
e-quasi compact space, e-compact spaces, sets. First, we introduce
the notion of e-open cover.

\begin{definition}
Let $(\widetilde{X},\tau,A)$ be a soft e-topological space,
$(F,A)\in S(\widetilde{X})$, and $\{(O_i,A)\}_{i\in I}$ be a family
 of soft e-open sets in $(\widetilde{X},\tau,A)$.
\begin{description}
    \item[i)] $\{(O_i,A)\}_{i\in I}$ is called a soft e-open
cover of $(\widetilde{X},A)$ if  $(\widetilde{X},A)= \underset{i\in
I}{\Cup}(O_i,A)$.
    \item[ii)]
$\{(O_i,A)\}_{i\in I}$ is called a soft e-open cover of $(F,A)$ if:
$(F,A)\widetilde{\subseteq} \underset{i\in I}{\Cup}(O_i,A)$.
\end{description}
\end{definition}

In the following definition, we introduce the notion of a soft
e-quasi compact space.

\begin{definition}
Let $(\widetilde{X},\tau,A)$ be a soft e-topological space.
$(\widetilde{X},\tau,A)$ is called a soft e-quasi compact space if
every soft e-open cover of $(\widetilde{X},A)$ has a finite
sub-e-cover of $(\widetilde{X},A)$.
\end{definition}

In the following theorem, we introduce a necessary condition, so
that a soft e-topological space has a soft e-quasi compact space.

\begin{theorem}\label{th4.1}
Let $(\widetilde{X},\tau,A)$ be a soft e-quasi compact; then for
every family $\{(F_i,A)\}_{i\in I}$ of soft e-closes such that
$\underset{i\in I}{\Cap}(F_i,A)=(\widetilde{\Phi},A),$ we can
extract a finite subfamily $\{(F_i,A)\}_{i\in I_0\subset I}$ such
that $\underset{i\in I_0}{\Cap}(F_i,A)=(\widetilde{\Phi},A).$
\end{theorem}

\noindent\emph{Proof.} Assume that $(\widetilde{X},\tau,A)$ is soft
e-quasi compact, and let $\{(F_i,A)\}_{i\in I}$ be a family of soft
e-closed such that $\underset{i\in
I}{\Cap}(F_i,A)=(\widetilde{\Phi},A).$ Hence,
$\{(F_i,A)^\mathbb{C}\}_{i\in I}$ is a family of soft e-opens, and
we have $\underset{i\in I}{\Cup}(F_i,A)^\mathbb{C}=
(\widetilde{X},A)$. Since $(\widetilde{X},\tau,A)$ is a quasi
e-compact, there exists $I_0\subset I$ such that $\underset{i\in
I_0}{\Cup}(F_i,A)\mathbb{^C}= (\widetilde{X},A)$. Then,
$\underset{i\in I_0}{\Cap}(F_i,A)= (\widetilde{\Phi},A).$

\begin{remark}
Since the complementary of a soft e-open set is not a soft e-close
set in general, the converse of theorem \ref{th4.1} is not true in
general. This is shown in the following counter example. Let
 $X=]1,+\infty[, A=[1,+\infty[, I=[1,+\infty[,$ and
$\tau=\{(\widetilde{\Phi},A)\}\cup\{(O_i,A), i\in I\}$, where:
$(O_i,A)(\alpha)=\displaystyle\left]\frac{1+i\alpha}{i+\alpha},+\infty\right[,$
for all $\alpha\in A$. First, it is obvious that
$(\widetilde{\Phi},A)\in \tau,$ and $(\widetilde{X},A)=(O_1,A)\in
\tau.$ On other hand, for any collection $I_0\subset I$ we have:
$\underset{i\in I_0}{\Cup}(O_i,A)= (O_{i_0},A),$ where
$i_0=\displaystyle\min\{i, i\in I_0\}$. Finally, for all $i,j\in I$
such that $i<j,$ we have: $(O_i,A)\Cap(O_j,A)=(O_j,A)$. Then,
$(\widetilde{X},\tau,A)$ is a soft e-topological space.

Now;  for all $i\in I: (O_i,A)(1)=X,$ then $O_i^c(1)=\Phi,$ and
$O_i^c(\alpha)\neq\Phi$ for all $\alpha\neq1$ i.e. $(O_i,A)^c\notin
S(\widetilde{X})$ for all $i\in ]1,+\infty[$. The collection of soft
e-closed is only $\{(\widetilde{\Phi},A), (\widetilde{X},A)\},$
which is finite,  but the family $\{(O_i,A)\}_{i\in ]1,+\infty[}$ is
a soft e-open cover of $(\widetilde{X},A)$, and we can't extract a
finite e-open subcover of $(\widetilde{X},A)$.
\end{remark}

In the following definition, we introduce the notion of a soft
e-compact space.

\begin{definition}
Let $(\widetilde{X},\tau,A)$ be a soft e- topological space.
$(\widetilde{X},\tau,A)$ is called a soft e-compact space, if
$(\widetilde{X},\tau,A)$ is soft quasi e-compact space, and soft
e-Hausdorff space.
\end{definition}

Now, we introduce a necessary condition, so that a soft
e-topological space has a soft e-compact space.

\begin{theorem}\label{th4.2}
Let $(\widetilde{X},\tau,A)$ be a soft e-compact space, and
$\{(F_i,A)\}_{i=1}^\infty$ be a family of decreasing soft e-closed
sets, then $\Cap_{i=1}^\infty(F_i,A)\neq(\widetilde{\Phi},A),$
\end{theorem}

\noindent\emph{Proof.} Assume that
$\Cap_{i=1}^\infty(F_i,A)=(\widetilde{\Phi},A),$ then
$\Cup_{i=1}^\infty(F_i,A)^{\mathbb{C}}=(\widetilde{X},A).$ Since
$(\widetilde{X},A)$ is a soft e-compact, and
$\{(F_i,A)^{\mathbb{C}}\}_{i=1}^\infty$ is a family of soft e-open
sets, we can extract a decreasing finite subfamily
$\{(F_{i_k},A)^{\mathbb{C}}\}_{k=1}^n$ such that
$\Cup_{k=1}^n(F_{i_k},A)^{\mathbb{C}}=(\widetilde{X},A).$ Then
$(F_{i_n},A)= \Cap_{k=1}^n
(F_{i_k},A)^{\mathbb{C}}=(\widetilde{\Phi},A),$ which is a
contradiction.

In the following definition, we introduce the notion of a soft
e-compact set.

\begin{definition}
Let $(\widetilde{X},\tau,A)$ be a soft e-Hausdorff space, and
$(F,A)\in S(\widetilde{X})$ such that $(F,A)^C\in S(\widetilde{X})$.
$(F,A)$ is called a soft e-compact set if all soft e-open cover of
$(F,A)$ has a finite sub-e-cover of $(F,A)$.
\end{definition}

\begin{example}\label{ex3.1}
Let $X=\mathbb{R}, A=\{\alpha,\beta\},$ and $\tau$ be the collection
of soft sets $(O,A)\in S(\widetilde{X})$ such that
$(O,A)=(\widetilde{\Phi},A),$ or for all $\widetilde{x}
\widetilde{\in} (F,A),$ there exist $r_\alpha>0, r_\beta>0$ such
that $]\widetilde{x}(\alpha)-r_\alpha, \widetilde{x}(\alpha)+
r_\alpha[\subset O(\alpha),$ and $]\widetilde{x}(\beta)-r_\beta,
\widetilde{x}(\beta)+r_\beta[\subset O(\beta)$. It is obvious that
$(\widetilde{X},\tau,A)$ is a soft e-topological space, and
$\tau_\alpha=\{O(\alpha), (O,A)\in \tau\}, \tau_\beta=\{O(\beta),
(O,A)\in \tau\}$ are two crisp topologies  of $X,$ equivalent to the
topology of metric space $(\mathbb{R},|.|)$. Now, let
$\widetilde{x}, \widetilde{y}$ be two soft elements of
$(\widetilde{X},A)$ such that $\widetilde{x}(\alpha)\neq
\widetilde{y}(\alpha), \widetilde{x}(\beta)\neq
\widetilde{y}(\beta),$ and $r_\alpha=|\widetilde{x}(\alpha)-
\widetilde{y}(\alpha)|, r_\beta=|\widetilde{x}(\beta)-
\widetilde{y}(\beta)|$ then $\left]\widetilde{x}(\alpha)-
\frac{r_\alpha}{3}, \widetilde{x}(\alpha)+
\frac{r_\alpha}{3}\right[\cap \left]\widetilde{y}(\alpha)-
\frac{r_\alpha}{3}, \widetilde{y}(\alpha)+
\frac{r_\alpha}{3}\right[=\emptyset,$ and
$\left]\widetilde{x}(\beta)- \frac{r_\beta}{3},
\widetilde{x}(\beta)+ \frac{r_\beta}{3}\right[\cap
\left]\widetilde{y}(\beta)- \frac{r_\beta}{3}, \widetilde{y}(\beta)+
\frac{r_\beta}{3}\right[=\emptyset.$ Hence, there exists $(F,A)$ a
soft e-nbd of $\widetilde{x},$ $(G,A)$ a soft e-nbd of
$\widetilde{y}$ such that
$(F,A)\widetilde{\cap}(F,A)=(\widetilde{\Phi},A)$\\
$(F(\alpha)=\left]\widetilde{x}(\alpha)- \frac{r_\alpha}{3},
\widetilde{x}(\alpha)+ \frac{r_\alpha}{3}\right[,
F(\beta)=\left]\widetilde{x}(\beta)- \frac{r_\beta}{3},
\widetilde{x}(\beta)+ \frac{r_\beta}{3}\right[,$\\
$G(\alpha)=\left]\widetilde{y}(\alpha)- \frac{r_\alpha}{3},
\widetilde{y}(\alpha)- \frac{r_\alpha}{3}\right[, G(\beta)=
\left]\widetilde{y}(\beta)- \frac{r_\beta}{3}, \widetilde{y}(\beta)-
\frac{r_\beta}{3}\right[$. Hence, $(\widetilde{X},\tau,A)$ is a soft
e-Hausdorff space. Finally, $(\widetilde{X},\tau,A)$ is not a soft
e-compact space, since we can't extract a soft e-sub-cover of the
e-cover open $\{(F_n,A), n\in \mathbb{N})\}
(F_n(\alpha)=F_n(\beta)=]-n,n[),$ but the soft set
$(\widetilde{[-1,1]},A)$ is a soft e-compact set.
\end{example}

In the following theorem, we show that the sub-e-compact space is a
soft e-compact set.

\begin{theorem}\label{th4.3}
Let $(\widetilde{X},\tau,A)$ be  a soft e-Hausdorff space. Assume
that for all $(O_1,A), (O_2,A)\in \tau$ we have:
$(O_1,A)\widetilde{\cap} (O_2,A)\in S(\widetilde{X})$. Let $Y$ be a
nonempty subset of $X$ such that $Y\neq X,$ and for all $(O,A)\in
\tau$ we have: $(O,A)\widetilde{\cap}(\widetilde{Y},A)\in
S(\widetilde{X})$. If $(\widetilde{Y},\tau_Y,A)$ is a soft e-compact
space, then $(\widetilde{Y},A)$ is a soft e-compact set in
$(\widetilde{X},\tau,A)$.
\end{theorem}

\noindent\emph{Proof.} Since $Y\neq\emptyset, Y\neq X,$ we obtain
$(\widetilde{Y},A), (\widetilde{Y},A)^C\in S(\widetilde{X})$. Now,
let $\{(O_i,A), i\in I\}\subset \tau$ such that $(\widetilde{Y},A)
\widetilde{\subseteq}\underset{i\in I}{\Cup}(O_i,A)$. Then
$(\widetilde{Y},A) =(\widetilde{Y},A)\Cap[\underset{i\in
I}{\Cup}(O_i,A)]= \underset{i\in
I}{\Cup}[(\widetilde{Y},A)\Cap(O_i,A)]$. Hence, the family
$\{(\widetilde{Y},A)\Cap(O_i,A),i\in I\}$ is a soft $e_Y$-open cover
of $(\widetilde{Y},A),$ and since $(\widetilde{Y},\tau_Y,A)$ is a
soft e-compact space we can extract a finite family
$\{(\widetilde{Y},A)\Cap(O_i,A),i\in I_0\}$ such that
$(\widetilde{Y},A)= \underset{i\in
I_0}{\Cup}[(\widetilde{Y},A)\Cap(O_i,A)]=
(\widetilde{Y},A)\Cap[\underset{i\in I_0}{\Cup}(O_i,A)]$. Then,
$(\widetilde{Y},A) \widetilde{\subseteq}\Cup_{i\in I_0}(O_i,A),$
hence $(\widetilde{Y},A)$ is a soft e-compact set in
$(\widetilde{X},\tau,A)$.

\begin{remark}
The converse of theorem \ref{th4.3} is not true in general. This is
shown in the following counter example. Consider the soft
e-topological space, which is introduced in example \ref{ex3.1}. Let
$X=\{a,b,c,d\}, A=\{\alpha,\beta\},$ and
$\tau=\{(\widetilde{\Phi},A), (F,A), (G,A), (H,A),
(\widetilde{X},A)\},$ such that $F(\alpha)=\{a\}, F(\beta)=\{c,d\},
G(\alpha)=\{c,d\}, G(\beta)=\{a\}, H(\alpha)=\{a,c,d\},
H(\beta)=\{a,c,d\}$. Let $Y=\{b,c,d\}$. $(\widetilde{Y},A)$ is a
soft compact set, since it is a finite soft set, but we can't
introduce a soft sub-e-topological space from  $(\widetilde{Y},A)$,
since $(F,A)\widetilde{\cap}(\widetilde{Y},A) \notin
S(\widetilde{X}).$
\end{remark}

In the following two theorems, we show  the relationship between
soft e-compact set and soft e-closed set.
\begin{theorem}\label{th4.4}
Let $(\widetilde{X},\tau,A)$ be  a soft e-Hausdorff  such that for
all $(O_1,A), (O_2,A)\in \tau$ we have: $(O_1,A)\widetilde{\cap}
(O_2,A)\in S(\widetilde{X})$. Let $(F,A)$ be a soft e-compact set,
then $(F,A)$ is a soft e-closed set.
\end{theorem}

\noindent\emph{Proof.} Assume that $(F,A)^{\mathbb{C}}\neq
(\widetilde{\Phi},A),$ and let $\widetilde{y}\widetilde{\in}
(F,A)^{\mathbb{C}},$ then for all $\widetilde{x}\widetilde{\in}
(F,A)$ we have $\widetilde{x}(\alpha)\neq\widetilde{y}(\alpha)$ for
all $\alpha\in A.$ Since $(\widetilde{X},\tau,A)$ is a soft
e-Hausdorff, there exist $(G_x,A), (H_x,A)\in \tau$ such that
$\widetilde{x}\widetilde{\in} (G_x,A), \widetilde{y}\widetilde{\in}
(H_x,A)$ and $(G_x,A)\widetilde{\cap}(H_x,A)= (\widetilde{\Phi},A).$
We have $(F,A)\widetilde{\subseteq}
\underset{\widetilde{x}\widetilde{\in} (F,A)}{\Cup} (G_x,A),$ and
since $(F,A)$ be a soft e-compact set, there exist $\widetilde{x}_1,
\widetilde{x}_2, \ldots, \widetilde{x}_n\widetilde{\in}(F,A)$ such
that $F,A)\widetilde{\subseteq} \Cup_{i=1}^n
(G_{x_i},A)=\Cup_{i=1}^n (G_i,A)$. Putting $(H,A)=\Cup_{i=1}^n
(H_{x_i},A)\in\tau,$ then $\widetilde{y}\in (H,A)$ and
$(H,A)\Cap(G_i,A)=(\widetilde{\Phi},A)$ for all $i=1\ldots n.$ Since
$(H,A), (G_i,A)\in \tau$, we have $(H,A)\Cap
(G_i,A)=(H,A)\widetilde{\cap} (G_i,A)$ for all $i=1\ldots n$. Then
for all $\alpha\in A$ we have $\Cup_{i=1}^n
[(H,A)\Cap(G_i,A)](\alpha)=\widetilde{\cup}_{i=1}^n
[(H,A)\widetilde{\cap}(G_i,A)](\alpha)=(H,A)\widetilde{\cap}
[\widetilde{\cup}_{i=1}^n (G_i,A)](\alpha)=\emptyset.$ Since $(H,A),
\widetilde{\cup}_{i=1}^n (G_i,A)\in \tau,$ then
$(H,A)\widetilde{\cap} [\widetilde{\cup}_{i=1}^n (G_i,A)]\in
S(\widetilde{X})$. Hence $(H,A)\widetilde{\cap}
[\widetilde{\cup}_{i=1}^n (G_i,A)]= (H,A)\Cap [\Cup_{i=1}^n
(G_i,A)]= (\widetilde{\Phi},A)$. Then;\\
 $(H,A)\widetilde{\subseteq}
(\Cup_{i=1}^n (G_i,A))^C\widetilde{\subseteq} (F,A)^C$, so
$(F,A)^\mathbb{C}\in \tau,$ and $(F,A)$ is soft e-closed set.

\begin{theorem}\label{th4.5}
Let $(\widetilde{X},\tau,A)$ be  a soft e-Hausdorff. Assume that
there exists a soft e-compact set $(K,A)$ such that
$(F,A)\widetilde{\subseteq} (K,A),$  and let $(F,A)$ be a soft
e-closed set, then $(F,A)$ is a soft e-compact set.
\end{theorem}

\noindent\emph{Proof.} Since $(F,A)$ is a soft e-closed set,
$(F,A)^C\in S(\widetilde{X})$ and $(F,A)^\mathbb{C}\in \tau.$ Let
$\{(O_i,A), i\in I\}$ be a soft e-cover open of $(F,A),$ then
$(\widetilde{X},A)= (F,A)^\mathbb{C}\Cup [\underset{i\in I}{\Cup}
(O_i,A)]$, and $(K,A)\widetilde{\subseteq} (F,A)^\mathbb{C}\Cup
[\underset{i\in I}{\Cup} (O_i,A)]$. Since $(K,A)$ is soft e-compact
set, there exists a finite subfamily $\{(O_i,A), i\in I_0\subset
I\}$ such that $(K,A)\widetilde{\subseteq} (F,A)^\mathbb{C}\Cup
[\underset{i\in I_0}{\Cup} (O_i,A)]$. Then,
$(F,A)\widetilde{\subseteq} \underset{i\in I_0}{\Cup} (O_i,A),$
hence $(F,A)$ is a soft e-compact set.

The following proposition provide that the soft e-compactness is
compatibly as the soft elementary union, and elementary intersection
(in more conditions).

\begin{proposition}\label{prop4.1}
Let $(\widetilde{X},\tau,A)$ be  a soft e-Hausdorff. Then
\begin{enumerate}
    \item elementary union of two soft e-compact sets is a soft
    e-compact set.
    \item if $(\widetilde{X},\tau,A)$ is a soft e-compact, and for all $(O_1,A),
(O_2,A)\in \tau$ we have: $(O_1,A)\widetilde{\cap} (O_2,A)\in
S(\widetilde{X})$, then elementary intersection of any soft
e-compact sets is a soft e-compact set.
\end{enumerate}
\end{proposition}

\noindent\emph{Proof.} Let $(\widetilde{X},\tau,A)$ be  a soft
e-Hausdorff space.
\begin{enumerate}
    \item Let $(K_1,A), (K_2,A)$ be two soft e-compact sets, and let $\{(O_i,A),i\in
    I\}$ be an e-open cover of $(K_1,A)\Cup(K_2,A),$ Then $\{(O_i,A),i\in
    I\}$ be an e-open cover of  $(K_1,A)$ and $(K_2,A).$
    We can extract a finite subcover $\{(O_i,A),i\in
    I_1\}$ of $(K_1,A)$ and a finite subcover $\{(O_i,A),i\in
    I_2\}$ of $(K_2,A)$.
    Hence; $\{(O_i,A),i\in I_1\cup I_2\}$ is a subcover of  $(K_1,A)\Cup(K_2,A).$
    \item Assume that $(\widetilde{X},\tau,A)$ is a soft e-compact and for all $(O_1,A),
(O_2,A)\in \tau$ we have: $(O_1,A)\widetilde{\cap} (O_2,A)\in
S(\widetilde{X})$. Let $\{(K_i,A),i\in I\}$ be a family of soft
e-compact sets, then $\{(K_i,A),i\in I\}$ be a family of soft
e-closed sets. $\underset{i\in I}{\Cap} (K_i,A)$ is a soft e-closed
set and soft subset of any soft e-compact set $(K_i,A),$ hence
$\underset{i\in I}{\Cap} (K_i,A)$ is a soft e-compact set.
\end{enumerate}

\begin{remark}\label{rem4.2}
We can replace the condition $(\widetilde{X},\tau,A)$ be a soft
e-compact by the condition: there exists a soft e-compact set
$(K,A)$ such that for all $i\in I$ we have
$(K_i,A)\widetilde{\subseteq}(K,A).$
\end{remark}

Now, we prove some properties of soft e-compact sets, and spaces.

\begin{theorem}\label{th4.6}
Let $(\widetilde{X},\tau,A)$ be  a soft e-Hausdorff, $(K,A)$ be a
soft e-compact and $(F,A)\in S(\widetilde{X})$ be a soft subset not
finite of $(K,A)$.  Then $(F,A)$ has a limiting soft element.
\end{theorem}

\noindent\emph{Proof.} Let $(\widetilde{X},\tau,A)$ be  a soft
e-Hausdorff,  $(K,A)$ be a soft e-compact, $(F,A)\in
S(\widetilde{X}).$ Assume that $(F,A)\widetilde{\subseteq}(K,A)$ is
not finite, and has not a soft limiting element, than by lemma
\ref{lem3.1} for all $\widetilde{x}\widetilde{\in}
(\widetilde{X},A),$ there exists $(G_x,A)\in \tau$ such that for all
$\widetilde{y}\widetilde{\in} (\widetilde{X},A),
\widetilde{y}\widetilde{\in} (F,A)\Cap(G_x,A)$ implies that
$\widetilde{y}=\widetilde{x}.$ The family $\{(G_x,A),
\widetilde{x}\widetilde{\in} (\widetilde{X},A)\}$ is a soft e-open
cover of $(K,A).$ So, we can extract a finite cover
$\{(G_{x_i},A),\widetilde{x_i}\widetilde{\in} (\widetilde{X},A),
i=1\ldots n\}= \{(G_i,A),\widetilde{x_i}\widetilde{\in}
(\widetilde{X},A), i=1\ldots n\}.$ In the family
$\{\widetilde{x_i}\widetilde{\in} (\widetilde{X},A), i=1\ldots n\}$
there exists a subfamily $\{\widetilde{x_i}\widetilde{\in}
(\widetilde{X},A), i\in \{1,2,\ldots, n\}\}$. Since,
$(K,A)\widetilde{\subseteq}
[\Cup_{i=1}^n(G_i,A)]_{\widetilde{x}_i\widetilde{\in}(F,A)} \Cup
[\Cup_{i=1}^n(G_i,A)]_{\widetilde{x}_i\widetilde{\notin}(F,A)}$.
Then,  $(F,A)$ is a soft subset of $SS(\{\widetilde{x}_i
\widetilde{\in}(F,A), i=1\ldots n\})$, hence $(F,A)$ is finite,
which is a contradiction.

\begin{theorem}\label{th4.7}
Let $(\widetilde{X},\tau,A)$ be  a soft e-Hausdorff. If
$(\widetilde{X},\tau,A)$ is a soft e-compact, then
$(\widetilde{X},\tau,A)$ is a soft e-regular space.
\end{theorem}

\noindent\emph{Proof.} Let $(F,A)$ be a soft e-closed  and let
$\widetilde{y}\widetilde{\in} (F,A)^\mathbb{C}$. From the proof of
theorem \ref{th4.4}, there exist $(G,A), (H,A)\in \tau$ such that
$(F,A)\widetilde{\subseteq} (G,A), \widetilde{y}\widetilde{\in}
(H,A),$ and $(G,A)\Cap(H,A)=(\widetilde{\Phi},A).$ Then;
$(\widetilde{X},\tau,A)$ is a soft e-regular space.

\begin{theorem}\label{th4.8}
Let $(\widetilde{X},\tau,A)$ be  a soft e-Hausdorff. If
$(\widetilde{X},\tau,A)$ is a soft e-compact, then
$(\widetilde{X},\tau,A)$ is a soft e-normal space.
\end{theorem}

\noindent\emph{Proof.} Let $(F,A)$ be a soft e-closed set, and let
$\widetilde{y}\widetilde{\in} (F,A)^\mathbb{C}$. From the proof of
theorem \ref{th4.6}, if $(F_1,A), (F_2,A)$ are two soft e-closed
sets such that $(F_1,A)\widetilde{\cap}(F_2,A)=
(\widetilde{\Phi},A)$, and for all $\widetilde{y}\widetilde{\in}
(F_2,A)$ there exist $(G_1^y,A), (G_2^y,A)\in \tau$ such that
$(F,A)\widetilde{\subseteq} (G_1^y,A), \widetilde{y}\widetilde{\in}
(G_2^y,A),$ and $(G_1^y,A)\Cap(G_2^y,A)=(\widetilde{\Phi},A).$ We
can extract a sub-e-cover $\{(G_2^{y_i},A)\}_{i=1}^n$ from the cover
$\{(G_2^y,A), \widetilde{y}\widetilde{\in} (F_2,A)\}$ of $(F_2,A)$.
Putting, $(G_1^i,A)= (G_1^{y_i},A), (G_2^i,A)= (G_2^{y_i},A)$ for
all $i=1\ldots n, (G_1,A)=\Cap_{i=1}^n(G_1^i,A),
(G_2,A)=\Cap_{i=1}^n(G_2^i,A),$ hence $(F_1,A)\widetilde{\subseteq}
(G_1,A), (F_2,A)\widetilde{\subseteq} (G_2,A)$ and
$(G_1,A)\Cap(G_2,A)= (\widetilde{\Phi},A)$. Then
$(\widetilde{X},\tau,A)$ is a soft e-normal space.

This proposition provide that the image of a soft e-compact set by a
soft e-continuous function in a soft e-Hausdorff space is a soft
e-compact set.

\begin{proposition}\label{prop6.1}
Let $(\widetilde{X},\tau,A), (\widetilde{Y},\sigma,A)$ be two soft
e-Hausdorff spaces. Let  $f: SE(\widetilde{X}) \rightarrow
SE(\widetilde{Y})$ be a soft function and $(K,A)$ be a soft
e-compact set of $(\widetilde{X},\tau,A)$. If $f$ is a soft
e-continuous then $f[(K,A)]$ is a soft e-compact set of
$(\widetilde{Y},\sigma,A)$.
\end{proposition}

\noindent\emph{Proof.} Let $\{(U_i,A)\in \sigma,i\in I\}$ be an open
cover of $f[(K,A)],$ then $f[(K,A)]\widetilde{\subseteq} \Cup_{i\in
I} (U_i,A).$ By Proposition 5.4 of \cite{[Chiney-Samanta]}, we have
$(K,A)\widetilde{\subseteq} f^{-1}(f[(K,A)])\widetilde{\subseteq}
f^{-1}[\underset{i\in I}{\Cup} (U_i,A)]=\underset{i\in
I}{\Cup}f^{-1} [(U_i,A)].$ Since $f$ is a soft e-continuous
function, then $\{f^{-1} [(U_i,A)],i\in I\}$ is a soft e-open cover
of $(K,A),$ and since $(K,A)$ is a soft e-compact set we can extract
a finite subcover $\{f^{-1} [(U_i,A)],i\in I_0\}$ of $(K,A),$ i.e.
$(K,A)\widetilde{\subseteq} \underset{i\in I_0}{\Cup}f^{-1}
[(U_i,A)]$. Then $f[(K,A)]\widetilde{\subseteq} f(\underset{i\in
I_0}{\Cup}f^{-1} [(U_i,A)]) = \underset{i\in I_0}{\Cup}f(f^{-1}
[(U_i,A)]) \widetilde{\subseteq} \underset{i\in I_0}{\Cup}(U_i,A) $.
Hence $f[(K,A)]$ is a soft e-compact set of
$(\widetilde{Y},\sigma,A).$

\vs 0.5cm

\section{Soft e-locally compact space and soft e-Baire theorem }\label{5}
In this section, we introduce the notion of soft e-locally compact
space, and we introduce a soft elementary version of the Baire
theorem. First, we introduce some definitions.

\begin{definition}
Let $(\widetilde{X},\tau,A)$ be  a soft e-Hausdorff.
$(\widetilde{X},\tau,A)$ is called a soft e-locally compact space
if, for all $\widetilde{x}\widetilde{\in}(\widetilde{X},A)$ and for
all soft neighborhood $(N,A)$ of $\widetilde{x},$ there exists a
soft e-compact neighborhood $(K,A)$ of $\widetilde{x},$ such that
$(K,A) \widetilde{\subseteq}(N,A)$.
\end{definition}

\begin{definition}
Let $(\widetilde{X},\tau,A)$ be  a soft e-Hausdorff, and  $(F,A)\in
S(\widetilde{X})$ such that $(F,A)\neq(\widetilde{\Phi},A)$.
\begin{description}
    \item[i)] $(F,A)$ is called soft e-nowhere dense $($or soft rare set $)$ if
$\overset{0}{\overline{(F,A)}}=(\widetilde{\Phi},A)$
    \item[iii)] $(F,A)$ is called of first e-category $($or soft meager set$)$ if it is a countable union of e-nowhere dense subsets
of $(\widetilde{X},A)$.
    \item[ii)] $(F,A)$ is called of seconde e-category if it is not of first e-category.
\end{description}
\end{definition}

\begin{definition}
Let $(\widetilde{X},\tau,A)$ be  a soft e- Hausdorff. We say that
$(\widetilde{X},\tau,A)$ is a soft e-Baire space if for all
countable family of soft e-closed sets $\{(F_i,A)\}_{i=1}^{\infty}$
such that $\overset{0}{\overbrace{(F_i,A)}}=(\widetilde{\Phi},A)$ we
have $\overset{0}{\overbrace{\Cup_{i=1}^{\infty}(F_i,A)}}=
(\widetilde{\Phi},A).$
\end{definition}

In this theorem, we give the soft elementary version of Baire
theorem, using the soft elementary intersection..

\begin{theorem}\label{th5.1}
Let $(\widetilde{X},\tau,A)$ be  a soft locally compact space such
that for all $(O_1,A), (O_2,A)\in \tau$ we have:
$(O_1,A)\widetilde{\cap} (O_2,A)\in S(\widetilde{X})$, then
$(\widetilde{X},\tau,A)$  is a soft Baire space.
\end{theorem}

\noindent\emph{Proof.} Let $(\widetilde{X},\tau,A)$ be  a soft
e-locally compact space. Let $\{(F_i,A)\}_{i=1}^{\infty}$ be a
family of soft e-closed sets such that
$\overset{0}{\overbrace{(F_i,A)}}=(\widetilde{\Phi},A)$ for all
$i=1,2, \ldots$. We set $\Cup_{i=1}^{\infty}(F_i,A)=(F,A)$. To prove
that $\overset{0}{\overbrace{(F_i,A)}}= (\widetilde{\Phi},A)$
 it is enough to prove
that for all $(O,A)\in\tau$ we have $(O,A)\Cap(F,A)^C\neq
(\widetilde{\Phi},A)$. Since
$\overset{0}{\overbrace{(F_1,A)}}=(\widetilde{\Phi},A)$ we have
$(O,A)\widetilde{\nsubseteq}(F_1,A),$ hence $(O,A)\Cap(F,A)^C\neq
(\widetilde{\Phi},A)$.Since $(O,A),(F_1,A)^C\in\tau,
(O,A)\widetilde{\cap} (F_1,A)\in S(\widetilde{X})$, there exists a
soft element  $\widetilde{x}_1\widetilde{\in} (O,A)\Cap(F_1,A)^C,$
and since $(\widetilde{X},\tau,A)$ is a soft e-locally compact space
there exists a soft e-compact set $(K_1,A) \widetilde{\subseteq}
(O,A)\Cap(F,A)^C$ such that $\widetilde{x}_1\widetilde{\in}
\overset{0}{\overbrace{(K_1,A)}}$. Next,
 $\overset{0}{\overbrace{(F_2,A)}}=(\widetilde{\Phi},A)$ and $\overset{0}{\overbrace{(K_1,A)}}\neq (\widetilde{\Phi},A),$
 then there exists a soft element $\widetilde{x}_2$  and a soft e-compact set $(K_2,A)$ such that $\widetilde{x}_2\widetilde{\in}
(K_2,A)\widetilde{\subseteq} \overset{0}{\overbrace{(K_1,A)}}
\Cap(F_2,A)^C$. Proche to proche we construct a countable family of
e-closed sets $\{(K_i,A)\}_{i=1}^\infty$ such that
$(K_1,A)\widetilde{\supseteq} (K_2,A)\widetilde{\supseteq} \ldots$.
The family $\{(K_i,A)\}_{i=1}^\infty$ is a family of decreasing soft
subsets of the soft e-compact set $(K_1,A)$. So,
$\{(K_i,A)\}_{i=1}^\infty$ is a family of decreasing soft e-closed
subsets of the soft e-compact set $(K_1,A)$, then by theorem
\ref{th4.2} we have $\Cap_{i=1}^\infty(K_i,A)
\neq(\widetilde{\Phi},A)$. Let $\widetilde{x}\widetilde{\in}
\Cap_{i=1}^\infty(K_i,A) \neq(\widetilde{\Phi},A),$ then
$\widetilde{x}\widetilde{\in} (K_i,A)\Cap(F_i,A)^C$ for all $i=1,2,
\ldots$. Hence, $\widetilde{x}\widetilde{\notin}(F_i,A)$ for all
$i=1,2, \ldots$. Then, $\widetilde{x} \widetilde{\notin}(F,A),$
which need to $(O,A)\nsubseteq (F,A),$ hence $
\overset{0}{\overbrace{(K_1,A)}}= (\widetilde{\Phi},A).$

\vs 0.5cm

\section{Conclusion}\label{6}

In this paper,  basing on the approach of Chiney and Samanta
\cite{[Chiney-Samanta]}, we have introduced a definition of soft
elementary compact set, and space. We have investigated some
properties of the soft elementary compactness, and we have proved
the main result, which is the soft elementary version of Baire
theorem.

\vs 0.5cm

\end{document}